\documentclass[italian,UKenglish,11pt]{amsart}
\usepackage{etex}

\usepackage[utf8]{inputenc}

\usepackage[a4paper,lmargin=2.2cm,rmargin=2.2cm,tmargin=3.0cm,bmargin=3.0cm]{geometry}

\usepackage{multirow}

\usepackage{babel}

\usepackage{rotating,color,hyperref,appendix}
\usepackage{tocvsec2}
\usepackage{booktabs}
\usepackage{float}
\usepackage{tikz}
\usetikzlibrary{positioning,arrows,shapes,decorations.pathmorphing,shadows}

\usepackage{amsfonts}

\newcommand{\Spin}[1]{\ensuremath{\text{\upshape\rmfamily Spin}(#1)}}

\newcommand{\Sp}[1]{\mathrm{Sp}(#1)}
\newcommand{\SO}[1]{\ensuremath{\text{\upshape\rmfamily SO}(#1)}}
\newcommand{\SU}[1]{\mathrm{SU}(#1)}
\newcommand{\U}[1]{\mathrm{U}(#1)}

\newcommand{\I}{\mathcal{I}}

\newcommand{\RH}{R^\HH}
\newcommand{\LH}{L^\HH}
\newcommand{\RO}{R}

\newcommand{\liespin}[1]{\mathop{\mathfrak{spin}}(#1)}

\newcommand{\Id}{\text{Id}}

\newcommand{\End}[1]{\mathrm{End}(#1)}
\newcommand{\Hom}[1]{\mathrm{Hom}(#1)}

\newcommand{\Endplus}[1]{\mathrm{End}^+(#1)}
\newcommand{\Endminus}[1]{\mathrm{End}^-(#1)}

\newcommand{\CC}{\mathbb{C}}
\newcommand{\HH}{\mathbb{H}}
\newcommand{\RR}{\mathbb{R}}

\newcommand{\OO}{\mathbb{O}}

\newcommand{\EIII}{\mathrm{E\,III}}

\newcommand{\EVI}{\mathrm{E\,VI}}

\newcommand{\EVIII}{\mathrm{E\,VIII}}

\newcommand{\Esix}{\mathrm{E}_6}
\newcommand{\Eseven}{\mathrm{E}_7}
\newcommand{\Eeight}{\mathrm{E}_8}

\numberwithin{equation}{section}

\theoremstyle{plain}
\newtheorem{te}{Theorem}[section]
\newtheorem*{te*}{Theorem}
\newtheorem{pr}[te]{Proposition}

\newtheorem{lm}[te]{Lemma}

\theoremstyle{definition}
\newtheorem{de}[te]{Definition}

\theoremstyle{remark}
\newtheorem{re}[te]{Remark}


\begin{document}

\begin{otherlanguage}{italian}
\author{Paolo Piccinni}
\address{Sapienza Universit\`a di Roma\\ Dipartimento di Matematica, piazzale Aldo Moro 2, I-00185, Roma, Italy}
\email{piccinni@mat.uniroma1.it}
\end{otherlanguage}

\title[On some Grassmannians]{On some Grassmannians\\ carrying an even Clifford structure}

\thanks{The author was supported by the GNSAGA group of INdAM, by the PRIN Research Project 2015 "Variet\`a reali e complesse: geometria, topologia e analisi armonica", and by the Research Project "Polynomial identities and combinatorial methods in algebraic and geometric structures" of Sapienza Universit\`a di Roma}

\maketitle

\begin{abstract}
We give an explicit description of the non-flat parallel even Clifford structures of rank 8, 6, 5 on some real, complex and quaternionic Grassmannians, and discuss the r\^ole of the octonions in them, in particular for some low dimensional examples.

\end{abstract}

\vspace{0.5cm}

\section{Introduction}

\noindent It is known for a long time that the Grassmannians 
\[
Gr_{4}(\RR^{m+4})=\frac{\mathrm{SO}(m+4)}{\mathrm{SO}(m)\times \mathrm{SO}(4)} , \qquad Gr_{2}(\CC^{m+2})= \frac{\mathrm{SU}(m+2)}{\mathrm{S}(\mathrm{U}(m)\times \mathrm{U}(2))},
\]
are examples of positive quaternion K\"ahler manifolds, and together with the projective spaces $\HH P^m$ the only (known) series of manifolds carrying such a structure. Their local compatible almost hypercomplex structures come visibly from their elements, oriented real $4$-planes and complex $2$-planes, and the hypercomplex structure on the planes extends to local almost hypercomplex structures on both series of Grassmannians via the isomorphism of vector bundles $$T Gr \cong W \otimes W^\perp,$$ involving the tautological $W$  and its orthogonal complement $W^\perp$ in the ambient linear space.  

Accordingly, one can ask whether the dimensional analogy allows to define some sort of "octonionic structure" on 
\begin{equation}\label{octograss}
Gr_{8}(\RR^{m+8})=\frac{\mathrm{SO}(m+8)}{\mathrm{SO}(m)\times \mathrm{SO}(8)},  \quad Gr_{4}(\CC^{m+4})=\frac{\mathrm{SU}(m+4)}{\mathrm{S}(\mathrm{U}(m)\times \mathrm{U}(4))}, \quad  Gr_{2}(\HH^{m+2})=\frac{\mathrm{Sp}(m+2)}{\mathrm{Sp}(m)\times \mathrm{Sp}(2)},
\end{equation}
and it is not surprising that these are the three series of Grassmannians appearing in the classification Table 2 of \cite[page 955]{ms}, collecting manifolds that admit a parallel non-flat even Clifford structure.

This refers to the notion of \emph{even Clifford structure} on Riemannian manifolds $(M,g)$, introduced by A. Moroianu and U. Semmelmann in the same paper \cite{ms}, and defined as the datum of a (possibly only locally defined) real oriented Euclidean vector bundle $E$ over $M$, together with an algebra bundle \emph{Clifford morphism} 
\[
\varphi: \; \text{Cl}^0 E \rightarrow \End{TM}
\] 
from the even Clifford bundle $\text{Cl}^0 E$, mapping the sub-bundle $\Lambda^2 E$ into skew-symmetric endomorphisms. The rank $r$ of $E$ is said to be the \emph{rank of the even Clifford structure}. The even Clifford structure $E$ is called \emph{parallel non-flat} if there exists a non-flat metric connection $\nabla^E$ on $E$ such that $\varphi$ is connection preserving, i.e.
\[\varphi(\nabla^E_X \sigma) = \nabla^g_X \varphi (\sigma),
\]
for every tangent vector $X \in TM$ and section $\sigma$ of $\text{Cl}^0 E$, and where $\nabla^g$ is the Levi Civita connection. 

When $r=2,3$ an even Clifford structure is equivalent to an almost Hermitian and an almost quaternion Hermitian structure (respectively K\"ahler and quaternion K\"ahler if parallel). Other choices $r \geq 4$ correspond to further interesting geometries fitting into the notion, and in the parallel non-flat hypothesis the highest possible ranks $r=10,12,16$ are achieved only on the exceptional symmetric spaces of compact type 
\[
\EIII = \frac{\Esix}{\mathrm{Spin}(10)\cdot\mathrm{U}(1)}, \quad  \EVI=  \frac{\Eseven}{\mathrm{Spin}(12)\cdot\mathrm{Sp}(1)}, \quad \EVIII = \frac{\Eeight}{\mathrm{Spin}^+(16)},
\]
and on their non compact duals. 

Cf. \cite{mp}, \cite{ah}, \cite{had} for further developments of the notion of even Clifford structures, and \cite{pp3}, \cite{ppv}, \cite{p} for the study of $\EIII,\EVI,\EVIII$ as \emph{octonionic K\"ahler manifolds}.

The present paper is aimed to understand how much the algebra $\OO$ of octonions allows to describe the mentioned even Clifford structures on the Grassmannians of the three series \eqref{octograss}. 

Concerning the first of the three series, we will examine in particular when the dimension of the ambient space is even:
$$Gr_{8}(\RR^{2n+8}) = \frac{\mathrm{SO}(2n+8)}{\mathrm{SO}(2n)\times \mathrm{SO}(8)},$$
so that the Grassmannian is a spin manifold (cf. \cite{bh}, \cite{str}), and its even Clifford structure can be defined by the a global vector bundle. 

Some interesting "low-dimensional" cases are:
\[
Gr_{8}(\RR^{10}) , \qquad Gr_{8}(\RR^{12}), \qquad Gr_{8}^\perp(\RR^{16})
\]
(the latter the $\mathbb Z_2$ quotient by the orthogonal complement involution $\perp$ on oriented $8$-planes), which appear as totally geodesic half-dimensional "octonionic" sub-manifolds of $\EIII , \EVI ,\EVIII$, respectively.
In particular $Gr_{8}(\RR^{10}) \cong Gr_{2}(\RR^{10})$, isometric to the complex non singular quadric $Q_8 \subset \CC P^9$ and a Hermitian symmetric space, can be looked at as the "projective line" $(\CC \otimes \OO) P^1$ over the complex octonions. This fact is of some help in the study of the geometry and topology of $\EIII$, that in turn can be viewed as the projective plane $(\CC \otimes \OO) P^2$ over the complex octonions, cf. \cite{im}, \cite{pp3}. Next, $Gr_{8}(\RR^{12}) \cong Gr_{4}(\RR^{12})$, is one of the mentioned quaternion K\"ahler Wolf spaces, and similarly it can be viewed as the "projective line" $(\HH \otimes \OO) P^1$ over the quaternionic octonions. Finally, $Gr_{8}^\perp(\RR^{16})$ is in some sense the "projective line" $(\OO \otimes \OO) P^1$ over the octonionic octonions.  Cf. \cite{cn}, \cite{e}, \cite[pages 195, 198, 200]{BaeOct}.

Going to the second series of Grassmannians, it is again convenient to refer to an even dimensional ambient space:
$$Gr_{4}(\CC^{2n+4}) = \frac{\mathrm{SU}(2n+4)}{\mathrm{S}(\mathrm{U}(2n)\times \mathrm{U}(4))},$$
that insures the Grassmannian to be spin \cite{bh}. Here the lowest dimensional case is $Gr_{4}(\CC^{6})$, a projective non singular variety in $\CC P^{14}$ known as the third Severi variety, cf. \cite{pp3}.

The Grassmannians of the third series are all spin, and we will discuss how $\Spin{5} \cong \Sp{2}$ is defining their even Clifford structure. Some details on $Gr_{2}(\HH^{4})$, sharing some features of the classical Klein quadric of projective geometry, are also given.

\vspace{0.3cm}

\section{$\Spin{8}$ and some of its subgroups}\label{Spin8}

We collect in this Section some informations on the group $\Spin{8} \subset \SO{16}$ and on the structure it defines on 16-dimensional Riemannian manifolds. We will deal also with its two subgroups $\Spin{6}$ and $\Spin{5}$. General references for the present approach are the file \cite{br} by R. Bryant and the book \cite[pages 271-289]{ha} by R. Harvey. We will adopt here the notations of \cite{br}, but the choice of generators of $\Spin{8}$ as in \cite{ha}. This is  coherent with our already quoted previous work \cite{pp1}, \cite{pp3}, \cite{ppv}, \cite{p}. 
\vspace{0.2cm}

\subsection{$\Spin{8} \subset \SO{16}$}. For any $u \in \OO$, look at the linear map 
\begin{equation}\label{HarSpC}
m_u: \OO \oplus \OO \rightarrow \OO \oplus \OO, \qquad m_u: \left(
\begin{array}{c}
x \\
x'
\end{array}
\right)
\longrightarrow
\left(
\begin{array}{cc}
0& R_u \\ 
-R_{\overline u} & 0
\end{array}
\right)  
\left(
\begin{array}{c} 
x \\ 
x'
\end{array}
\right),
\end{equation}
where $R_u, R_{\overline u}$ denote the right multiplications by the octonions $u, \overline u$. Then $$(m_u)^2=-\vert u \vert^2 \Id,$$ so that $m_u$ induces a representation of the Clifford algebra $\text{Cl} \, \OO$, generated by $\OO$ with its standard quadratic form, on the vector space $\OO \oplus \OO \cong \RR^{16}$. Since $\text{Cl} \, \OO$ has no bilateral non trivial ideals, the representation is faithful and a count of dimensions then shows that $m_u$ yields the isomorphism $\text{Cl} \, \OO \cong M_{16}(\RR).$ 

The subgroup $\Spin{8} \subset \SO{16} \subset \text{Cl} \, \OO$ is generated by the compositions $m_u \, m_v$ with $u,v \in S^7 \subset \OO$, represented by matrices
\begin{equation}\label{muv}
m_{u,v}= \left(
\begin{array}{cc}
-R_{u} \circ R_{\overline v} & 0\\ 
0& -R_{\overline u} \circ R_v 
\end{array}
\right) .
\end{equation}

Recall that when $u,v$ are orthogonal octonions, for any $z \in \OO$ one has the identity $(z \overline v)u=-(z\overline u)v$, cf. \cite[Formula A.7.c, page 142]{hl}. By applying this to the matrix \eqref{muv}, one gets for orthogonal octonions $u,v$
\begin{equation}\label{anticommute}
m_{v,u} =-m_{u,v},
\end{equation}
so that for orthonormal $u,v$:
\begin{equation}\label{acs}
m^2_{u,v} =- \Id,
\end{equation}
i.e. \emph{any $m_u$ with unitary $u$ and any $m_{u,v}: \OO \oplus \OO \rightarrow \OO \oplus \OO$ with orthonormal $u,v$ is a complex structure.}

Look now at the choices of $u$ in the standard basis of $\OO$:
\[
u = 1,i,j,k,e,f,g,h \in S^7 \subset \OO.
\]
The corresponding $m_u$ give the complex structures in $\RR^{16}$:
\begin{equation}\label{top}
\scriptsize{
\begin{aligned}
&m_1=\left(
\begin{array}{cc}
0 & \Id\\
-\Id & 0
\end{array}
\right), \quad
m_i=\left(
\begin{array}{cc}
0 & R_i\\
R_i & 0
\end{array}
\right), \quad
m_j=\left(
\begin{array}{cc}
0 & R_j\\
R_j& 0
\end{array}
\right), \quad
m_k=\left(
\begin{array}{cc}
0 & R_k\\
R_k & 0
\end{array}
\right),&
\\
&m_e=\left(
\begin{array}{cc}
0 & R_e\\
R_e & 0
\end{array}
\right), \quad
m_f=\left(
\begin{array}{cc}
0 & R_f\\
R_f& 0
\end{array}
\right), \quad 
m_g=\left(
\begin{array}{cc}
0 & R_g\\
R_g & 0
\end{array}
\right), \quad 
m_h=\left(
\begin{array}{cc}
0 & R_h\\
R_h & 0
\end{array}
\right),&
\end{aligned}
}
\end{equation}
where $R_i, R_j, \dots,R_h$ are the right multiplications by the unit octonions $i,j, \dots,h$. 

Their compositions $m_{u,v} = m_u m_v =-m_{v,u}$  $(u, v \in \{1,i,j,k,e,f,g,h\}, u \neq v)$, produce further 28 complex structures. 
Using notations $\RO_{\lambda, \mu} = \RO_\lambda \circ \RO_\mu$, for $\lambda,\mu\in\{i,j,k,e,f,g,h\}$, they are: 

\begin{equation}\label{eq:J1}
\scriptsize{
\begin{aligned}
m_{1,i}&=\left(
\begin{array}{cc}
\RO_i & 0 \\ 
0 & -\RO_i
\end{array}
\right),\, &
m_{1,j}&=\left(
\begin{array}{cc}
\RO_j & 0 \\ 
0 & -\RO_j
\end{array}
\right),\, &
m_{1,k}&=\left(
\begin{array}{cc}
\RO_k & 0 \\ 
0 & -\RO_k
\end{array}
\right), \, &
m_{1,e}&=\left(
\begin{array}{cc}
\RO_e & 0 \\ 
0 & -\RO_e
\end{array}
\right), 
\\
m_{1,f}&=\left(
\begin{array}{cc}
\RO_f & 0 \\ 
0 & -\RO_f
\end{array}
\right),\, &
m_{1,g}&=\left(
\begin{array}{cc}
\RO_g & 0 \\ 
0 & -\RO_g
\end{array}
\right),\, &
m_{1,h}&=\left(
\begin{array}{cc}
\RO_h & 0 \\ 
0 & -\RO_h
\end{array}
\right),\, &
m_{i,j}&=\left(
\begin{array}{cc}
\RO_{i,j} & 0 \\ 
0 & \RO_{i,j}
\end{array}
\right),
\\
m_{i,k}&=\left(
\begin{array}{cc}
\RO_{i,k} & 0 \\ 
0 & \RO_{i,k}
\end{array}
\right),\, &
m_{i,e}&=\left(
\begin{array}{cc}
\RO_{i,e} & 0 \\ 
0 & \RO_{i,e}
\end{array}
\right),\, &
m_{i,f}&=\left(
\begin{array}{cc}
\RO_{i,f} & 0 \\ 
0 & \RO_{i,f}
\end{array}
\right),\, &
m_{i,g}&=\left(
\begin{array}{cc}
\RO_{i,g} & 0 \\ 
0 & \RO_{i,g}
\end{array}
\right), 
\\
m_{i,h}&=\left(
\begin{array}{cc}
\RO_{i,h} & 0 \\ 
0 & \RO_{i,h}
\end{array}
\right),\, &
m_{j,k}&=\left(
\begin{array}{cc}
\RO_{j,k} & 0 \\ 
0 & \RO_{j,k}
\end{array}
\right),\, &
m_{j,e}&=\left(
\begin{array}{cc}
\RO_{j,e} & 0 \\ 
0 & \RO_{j,e}
\end{array}
\right),\, &
m_{j,f}&=\left(
\begin{array}{cc}
\RO_{j,f} & 0 \\ 
0 & \RO_{j,f}
\end{array}
\right), 
\\
m_{j,g}&=\left(
\begin{array}{cc}
\RO_{j,g} & 0 \\ 
0 & \RO_{j,g}
\end{array}
\right),\, &
m_{j,h}&=\left(
\begin{array}{cc}
\RO_{j,h} & 0 \\ 
0 & \RO_{j,h}
\end{array}
\right),\, &
m_{k,e}&=\left(
\begin{array}{cc}
\RO_{k,e} & 0 \\ 
0 & \RO_{k,e}
\end{array}
\right),\, &
m_{k,f}&=\left(
\begin{array}{cc}
\RO_{k,f} & 0 \\ 
0 & \RO_{k,f}
\end{array}
\right), 
\\
m_{k,g}&=\left(
\begin{array}{cc}
\RO_{k,g} & 0 \\ 
0 & \RO_{k,g}
\end{array}
\right),\, &
m_{k,h}&=\left(
\begin{array}{cc}
\RO_{k,h} & 0 \\ 
0 & \RO_{k,h}
\end{array}
\right),\, &
m_{e,f}&=\left(
\begin{array}{cc}
\RO_{e,f} & 0 \\ 
0 & \RO_{e,f}
\end{array}
\right),\, &
m_{e,g}&=\left(
\begin{array}{cc}
\RO_{e,g} & 0 \\ 
0 & \RO_{e,g}
\end{array}
\right),
\\
m_{e,h}&=\left(
\begin{array}{cc}
\RO_{e,h} & 0 \\ 
0 & \RO_{e,h}
\end{array}
\right), \, &
m_{f,g}&=\left(
\begin{array}{cc}
\RO_{f,g} & 0 \\ 
0 & \RO_{f,g}
\end{array}
\right),\, &
m_{f,h}&=\left(
\begin{array}{cc}
\RO_{f,h} & 0 \\ 
0 & \RO_{f,h}
\end{array}
\right),\, &
m_{g,h}&=\left(
\begin{array}{cc}
\RO_{g,h} & 0 \\ 
0 & \RO_{g,h}
\end{array}
\right),
\end{aligned}
}
\end{equation}
that are a basis of the Lie subalgebra $\frak{spin}(8) \subset \frak{so}(16)$. 

\begin{re}\label{triality}
Matrices in the subalgebra $\liespin{8} \subset \frak{so}(16)$ are characterized through the infinitesimal triality principle as block matrices:
\[
M=\left(
\begin{array}{cc}
m_+ & 0 \\ 
0& m_-
\end{array}
\right), 
\]
where $(m_+,m_-,m_0) \in \mathfrak{so}(8)$ are \emph{triality companions}, i.e. for each $u \in \OO$ there exists \mbox{$v=m_0(u)$} such that $R_v + R_u \, m_-= m_+ \, R_u$, cf. \cite[pages 278--279, 285]{ha}, \cite[page 189]{mu}. Moreover, matrices in  $\mathfrak{spin}_\Delta(7) \subset \liespin{8}$ are characterized as those with $m_+=m_-$. Thus, the 21 matrices that in \eqref{eq:J1} have both indices in $\{i,j,k,e,f,g,h\}$ are a basis of $\mathfrak{spin}_\Delta(7)$, whereas the 7 further matrices, with first index $1$, complete this basis in $\liespin{8}$.
\end{re}

The matrices listed in \eqref{eq:J1} will be used in Section \ref{realgrassmannians} to define the even Clifford structure on Grassmannians $Gr_{8}(\RR^{2n+8})$, and it is worth to explicit them as real matrices. Denoting by $\RH_i,\RH_j,\RH_k$, $\LH_i,\LH_j,\LH_k$ the right and the left multiplication operators on quaternions, note first that:
\begin{equation}\label{matricesspin7}
\scriptsize{
\begin{aligned}
&\RO_i=\left(
\begin{array}{cc}
\RH_i & 0 \\ 
0 & -\RH_i
\end{array}
\right),\; \qquad \qquad 
\RO_j=\left(
\begin{array}{cc}
\RH_j & 0 \\ 
0 & -\RH_j
\end{array}
\right),\; \qquad \qquad
\RO_k=\left(
\begin{array}{cc}
\RH_k & 0 \\ 
0 & -\RH_k
\end{array}
\right), & 
\\
&\RO_e=\left(
\begin{array}{cc}
0 & -\Id \\ 
\Id & 0
\end{array}
\right),\quad 
\RO_f=\left(
\begin{array}{cc}
0 & \LH_i \\ 
\LH_i & 0
\end{array}
\right), \quad
\RO_g=\left(
\begin{array}{cc}
0 & \LH_j \\ 
\LH_j & 0
\end{array}
\right), \quad
\RO_h=\left(
\begin{array}{cc}
0 & \LH_k \\ 
\LH_k & 0
\end{array}
\right), &
\end{aligned}
}
\end{equation}
where
\begin{equation}\label{eq:right}
\tiny{
\RH_i=\left(
\begin{array}{rrrr}
0 & -1 & 0 & 0 \\
1 & 0 & 0 & 0 \\
0 & 0 & 0 & 1 \\
0 & 0 & -1 & 0
\end{array}
\right)\enspace,\qquad
\RH_j=\left(
\begin{array}{rrrr}
0 & 0 & -1 & 0 \\
0 & 0 & 0 & -1 \\
1 & 0 & 0 & 0 \\
0 & 1 & 0 & 0
\end{array}
\right)\enspace,\qquad
\RH_k=\left(
\begin{array}{rrrr}
0 & 0 & 0 & -1 \\
0 & 0 & 1 & 0 \\
0 & -1 & 0 & 0 \\
1 & 0 & 0 & 0
\end{array}
\right)\enspace,
}
\end{equation}
and \begin{equation}\label{eq:left}
\tiny{
\LH_i=\left(
\begin{array}{rrrr}
0 & -1 & 0 & 0 \\
1 & 0 & 0 & 0 \\
0 & 0 & 0 & -1 \\
0 & 0 & 1 & 0
\end{array}
\right)\enspace,\qquad
\LH_j=\left(
\begin{array}{rrrr}
0 & 0 & -1 & 0 \\
0 & 0 & 0 & 1 \\
1 & 0 & 0 & 0 \\
0 & -1 & 0 & 0
\end{array}
\right)\enspace,\qquad
\LH_k=\left(
\begin{array}{rrrr}
0 & 0 & 0 & -1 \\
0 & 0 & -1 & 0 \\
0 & 1 & 0 & 0 \\
1 & 0 & 0 & 0
\end{array}
\right).
}
\end{equation}

Next:
\begin{equation}\label{21}
\scriptsize{
\begin{aligned}
\RO_{i,j}&=\left(
\begin{array}{cc}
-\RH_k & 0 \\ 
0 & -\RH_k
\end{array}
\right)\enspace,\quad &
\RO_{i,k}&=\left(
\begin{array}{cc}
\RH_j & 0 \\ 
0 & \RH_j
\end{array}
\right)\enspace,\quad &
\RO_{i,e}&=\left(
\begin{array}{cc}
0 & -\RH_i \\ 
-\RH_i & 0
\end{array}
\right)\enspace, \\ 
\RO_{i,f}&=\left(
\begin{array}{cc}
0 & \RH_i\LH_i \\ 
-\RH_i\LH_i & 0
\end{array}
\right)\enspace,\quad &
\RO_{i,g}&=\left(
\begin{array}{cc}
0 & \RH_i\LH_j \\ 
-\RH_i\LH_j & 0
\end{array}
\right)\enspace,\quad &
\RO_{i,h}&=\left(
\begin{array}{cc}
0 & \RH_i\LH_k \\ 
-\RH_i\LH_k & 0
\end{array}
\right)\enspace, \\ 
\RO_{j,k}&=\left(
\begin{array}{cc}
-\RH_i & 0 \\ 
0 & -\RH_i
\end{array}
\right)\enspace,\quad &
\RO_{j,e}&=\left(
\begin{array}{cc}
0 & -\RH_j \\ 
-\RH_j & 0
\end{array}
\right)\enspace,\quad &
\RO_{j,f}&=\left(
\begin{array}{cc}
0 & \RH_j\LH_i \\ 
-\RH_j\LH_i & 0
\end{array}
\right)\enspace, \\
\RO_{j,g}&=\left(
\begin{array}{cc}
0 & \RH_j\LH_j \\ 
-\RH_j\LH_j & 0
\end{array}
\right)\enspace,\quad &
\RO_{j,h}&=\left(
\begin{array}{cc}
0 & \RH_j\LH_k \\ 
-\RH_j\LH_k & 0
\end{array}
\right)\enspace,\quad &
\RO_{k,e}&=\left(
\begin{array}{cc}
0 & -\RH_k \\ 
-\RH_k & 0
\end{array}
\right)\enspace, \\
\RO_{k,f}&=\left(
\begin{array}{cc}
0 & \RH_k\LH_i \\ 
-\RH_k\LH_i & 0
\end{array}
\right)\enspace,\quad &
\RO_{k,g}&=\left(
\begin{array}{cc}
0 & \RH_k\LH_j \\ 
-\RH_k\LH_j & 0
\end{array}
\right)\enspace,\quad &
\RO_{k,h}&=\left(
\begin{array}{cc}
0 & \RH_k\LH_k \\ 
-\RH_k\LH_k & 0
\end{array}
\right)\enspace, \\
\RO_{ef}&=\left(
\begin{array}{cc}
-\LH_i & 0 \\ 
0 & \LH_i
\end{array}
\right)\enspace,\quad &
\RO_{e,g}&=\left(
\begin{array}{cc}
-\LH_j & 0 \\ 
0 & \LH_j
\end{array}
\right)\enspace,\quad &
\RO_{e,h}&=\left(
\begin{array}{cc}
-\LH_k & 0 \\ 
0 & \LH_k
\end{array}
\right)\enspace, \\
\RO_{f,g}&=\left(
\begin{array}{cc}
\LH_k & 0 \\ 
0 & \LH_k
\end{array}
\right)\enspace,\quad &
\RO_{f,h}&=\left(
\begin{array}{cc}
-\LH_j & 0 \\ 
0 & -\LH_j
\end{array}
\right)\enspace,\quad &
\RO_{g,h}&=\left(
\begin{array}{cc}
\LH_i & 0 \\ 
0 & \LH_i
\end{array}
\right)\enspace, \\
\end{aligned}
}
\end{equation}
where the compositions $\RH \LH$ read:
\begin{equation}\label{RL}
\tiny{
\begin{aligned}
\RH_i\LH_i&=\left(
\begin{array}{rrrr}
-1 & 0 & 0 & 0 \\
0 &-1 & 0 & 0 \\
0 & 0 & 1 & 0 \\
0 & 0 & 0 & 1
\end{array}
\right)\thinspace,\enspace &
\RH_i\LH_j&=\left(
\begin{array}{rrrr}
0 & 0 & 0 & -1 \\
0 & 0 & -1 & 0 \\
0 &-1 & 0 & 0 \\
-1 & 0 & 0 & 0
\end{array}
\right)\thinspace,\enspace &
\RH_i\LH_k&=\left(
\begin{array}{rrrr}
0 & 0 & 1 & 0 \\
0 & 0 & 0 & -1 \\
1 & 0 & 0 & 0 \\
0 &-1 & 0 & 0
\end{array}
\right)\thinspace, \\
\RH_j\LH_i&=\left(
\begin{array}{rrrr}
0 & 0 & 0 & -1 \\
0 & 0 & 1 & 0 \\
0 & 1 & 0 & 0 \\
-1 & 0 & 0 & 0
\end{array}
\right)\thinspace,\enspace &
\RH_j\LH_j&=\left(
\begin{array}{rrrr}
1 & 0 & 0 & 0 \\
0 & -1 & 0 & 0 \\
0 & 0 & 1 & 0 \\
0 & 0 & 0 & -1
\end{array}
\right)\thinspace,\enspace &
\RH_j\LH_k&=\left(
\begin{array}{rrrr}
0 & 1 & 0 & 0 \\
1 & 0 & 0 & 0 \\
0 & 0 & 0 & 1 \\
0 & 0 & 1 & 0
\end{array}
\right)\thinspace, \\
\RH_k\LH_i&=\left(
\begin{array}{rrrr}
0 & 0 & -1 & 0 \\
0 & 0 & 0 & -1 \\
-1 & 0 & 0 & 0 \\
0 & -1 & 0 & 0
\end{array}
\right)\thinspace,\enspace &
\RH_k\LH_j&=\left(
\begin{array}{rrrr}
0 & 1 & 0 & 0 \\
1 & 0 & 0 & 0 \\
0 & 0 & 0 & -1 \\
0 & 0 & -1 & 0
\end{array}
\right)\thinspace,\enspace &
\RH_k\LH_k&=\left(
\begin{array}{rrrr}
-1 & 0 & 0 & 0 \\
0 & 1 & 0 & 0 \\
0 & 0 & 1 & 0 \\
0 & 0 & 0 & -1
\end{array}
\right)\thinspace.
\end{aligned}
}
\end{equation}

\begin{re} \label{re} The complex structures $m_u$ and $m_{u,v}$, with $u,v \in \{1,i,j,k,e,f,g,h\}$ can be easily compared with the ones associated with the standard $\Spin{9}$ structures in $\RR^{16}$. By referring to the notations $J_{\alpha \beta}$ used  in the paper \cite{pp1}, one sees that the present $m_1,m_i,m_j,m_k,m_e,m_f,m_g,m_h$ coincide respectively with $J_{19},J_{29},J_{39},J_{49},J_{59},J_{69},J_{79}, J_{89}$, and that the 28 generators $m_{u,v}$ of $\frak{spin}(8)$ of formulas \eqref{eq:J1} with the $J_{\alpha \beta}$, $1 \leq \alpha <\beta \leq 8$.
\end{re}

We mention also the following alternative approach to the construction of the basis \eqref{eq:J1} of $\frak{spin}(8)$, using a \emph{Clifford system}, i.e. a system of anti-commuting self-dual involutions. Real dimension 16 allows irreducible choices of such Clifford systems with 6,7,8,9 involutions, cf. the general discussion in \cite{ppv} and in its references.

\begin{pr}\label{8involutions} The self-dual anti-commuting involutions
\begin{equation}\label{eq:IO}
\scriptsize{ \begin{aligned}
\I_1&=\left(
\begin{array}{cc}
0 & -\mathrm{Id} \\ 
-\mathrm{Id} & 0
\end{array}
\right)\enspace,\, &
\I_i&=\left(
\begin{array}{cc}
0 & -\RO_i \\ 
\RO_i & 0
\end{array}
\right)\enspace,\,&
\I_j&=\left(
\begin{array}{cc}
0 & -\RO_j \\ 
\RO_j & 0
\end{array}
\right)\enspace, \, &
\I_k&=\left(
\begin{array}{cc}
0 & -\RO_k \\ 
\RO_k & 0
\end{array}
\right)\enspace,\\
\I_e&=\left(
\begin{array}{cc}
0 & -\RO_e \\ 
\RO_e & 0
\end{array}
\right)\enspace,\, &
\I_f&=\left(
\begin{array}{cc}
0 & -\RO_f\\ 
\RO_f & 0
\end{array}
\right)\enspace, \, &
\I_g&=\left(
\begin{array}{cc}
0 & -\RO_g \\ 
\RO_g & 0
\end{array}
\right)\enspace,\, &
\I_h&=\left(
\begin{array}{cc}
0 & -\RO_h \\ 
\RO_h & 0
\end{array}
\right)\enspace
\end{aligned}
}
\end{equation}
yield as compositions of pairs of them the basis of $\frak{spin}(8)$ given by the 28 complex structures $$m_{u,v}= - \I_u \circ \I_v$$ 
corresponding to choices of distinct and "increasing" indices $u,v \in \{1,i,j,k,e,f,g,h\}$.
\end{pr}

\vspace{0.2cm}
\subsection{The subgroup $\Spin{6} \subset \Spin{8}$} A way to describe subgroups of $\Spin{8}$ is to restrict to some of the complex structures \eqref{top}. For example, if we exclude just $m_1$, maintaining on $\RR^{16}$ the seven complex structures $m_i,m_j,m_k,m_e,m_f,m_g,m_h$, we see them generating the diagonal $\mathrm{Spin_\Delta}(7) \subset \Spin{8}$ and their compositions $m_{i,j},m_{i,k}, \dots , m_{g,h}$, explicitly given by 21 of the matrices \eqref{eq:J1}, generate the Lie subalgebra $\mathfrak{spin}_\Delta (7) \subset \mathfrak{so} (16)$.

We are also interested in the diagonal subgroup $\mathrm{Spin_\Delta}(6) \subset \Spin{8} \subset \SO{16}$. To describe the structure it defines, we can choose six of the above $m_u$ among the \eqref{top}, and accordingly, the corresponding 15 complex structures among the \eqref{eq:J1}, that generate the Lie algebra $\mathfrak{spin}_\Delta (6)$. A convenient choice is to select $m_1,m_i,m_j,m_k,m_e, m_f$, since by looking at the sequence of formulas \eqref{top}, \eqref{matricesspin7}, \eqref{eq:right}, \eqref{eq:left}, we see that these are the six of the \eqref{top} that commute with the standard complex structure of $\RR^{16}$, and they can in fact written as matrices in $\SU{4} \cong \Spin{6}$, thus acting on $\CC^8$. We will use them to describe the even Clifford structure on Grassmannians $Gr_4 (\CC^{2n+4})$. Explicitly, the first five of the \eqref{matricesspin7} as matrices in $\SU{4}$ read:
\begin{equation}\label{complex}
\tiny{\begin{aligned}
\RO_i=\left(
\begin{array}{rrrr}
i & 0&0&0 \\ 
0& -i&0&0 \\
0&0&-i&0\\
0&0&0&i
\end{array}
\right), \quad
\RO_j=\left(
\begin{array}{rrrr}
0 & -1 &0&0\\ 
1& 0&0&0\\
0&0&0&1\\
0&0&-1&0
\end{array}
\right), \quad
\RO_k=\left(
\begin{array}{rrrr}
0 & i&0&0 \\ 
i&0&0&0\\
0&0&0&-i\\
0&0&-i&0
\end{array}
\right), \\
\RO_e=\left(
\begin{array}{rrrr}
0 & 0&-1&0\\ 
0&0&0&-1\\
1&0&0&0\\
0&1&0&0
\end{array}
\right), \quad
\RO_f=\left(
\begin{array}{rrrr}
0 & 0&i&0\\ 
0&0&0&i\\
i&0&0&0\\
0&i&0&0
\end{array}
\right). \qquad \qquad \qquad \qquad
\end{aligned}
}
\end{equation}

The Clifford systems' approach gives here the following:

\begin{pr}\label{6involutions} The self-dual anti-commuting involutions
\begin{equation}\label{eq:IO}
\scriptsize{ \begin{aligned}
\I_1=\left(
\begin{array}{cc}
0 & -\mathrm{Id} \\ 
-\mathrm{Id} & 0
\end{array}
\right),\, \quad
\I_i=\left(
\begin{array}{cc}
0 & -\RO_i \\ 
\RO_i & 0
\end{array}
\right),\, \quad
\I_j=\left(
\begin{array}{cc}
0 & -\RO_j \\ 
\RO_j & 0
\end{array}
\right), \, \\
\I_k=\left(
\begin{array}{cc}
0 & -\RO_k \\ 
\RO_k & 0
\end{array}
\right), \, \quad
\I_e=\left(
\begin{array}{cc}
0 & -\RO_e \\ 
\RO_e & 0
\end{array}
\right),\, \quad
\I_f=\left(
\begin{array}{cc}
0 & -\RO_f\\ 
\RO_f & 0
\end{array}
\right)
\end{aligned}
}
\end{equation}
yield as compositions of pairs of them the basis of $\frak{spin}_\Delta(6)$ given by the 15 complex structures $$m_{u,v}= - \I_u \circ \I_v,$$ 
corresponding to choices of distinct and "increasing" indices $u,v \in \{1,i,j,k,e,f\}$.
\end{pr}

\vspace{0.2cm}
\subsection{$\Spin{5}$} It is convenient here to look at the representation $\Spin{5} \cong \Sp{2} \subset \SO{8}$. A coherent approach is via five anti-commuting self dual involutions of $\HH^2 \cong \RR^8$, whose compositions give rise to a basis of $\mathfrak{spin}(5) \cong \mathfrak{sp}(2)$. The simplest choices for such involutions are the following Pauli type matrices, cf \cite{pp1}:
\begin{equation}\label{eq:IH}
\scriptsize{\begin{aligned}
\sigma_1=\left(
\begin{array}{c|c}
0 & \Id \\ \hline
\Id & 0
\end{array}
\right), \quad
\sigma_2=\left(
\begin{array}{c|c}
0 & -\RH_i \\ \hline
\RH_i & 0
\end{array}
\right), \quad
\sigma_3=\left(
\begin{array}{c|c}
0 & -\RH_j \\ \hline
\RH_j & 0
\end{array}
\right), \quad
\sigma_4=\left(
\begin{array}{c|c}
0 & -\RH_k \\ \hline
\RH_k & 0
\end{array}
\right), \quad
\sigma_5=\left(
\begin{array}{c|c}
\Id & 0 \\ \hline
0 & -\Id
\end{array}
\right),
\end{aligned}
}
\end{equation}
whose compositions $\sigma_{\alpha\beta} = \sigma_\alpha\sigma_\beta$, $(\alpha<\beta)$ read
\begin{equation}\label{eq:Jspin51}
\scriptsize{\begin{aligned}
\sigma_{12}&=\left(
\begin{array}{c|c}
\RH_i & 0 \\ \hline
0 & -\RH_i
\end{array}
\right)\enspace,\qquad &
\sigma_{13}&=\left(	
\begin{array}{c|c}
\RH_j & 0 \\ \hline
0 & -\RH_j
\end{array}
\right)\enspace,\qquad &
\sigma_{14}&=\left(
\begin{array}{c|c}
\RH_k & 0 \\ \hline
0 & -\RH_k
\end{array}
\right)\enspace, \\
\sigma_{23}&=\left(
\begin{array}{c|c}
\RH_k & 0 \\ \hline
0 & \RH_k
\end{array}
\right)\enspace,\qquad &
\sigma_{24}&=\left(
\begin{array}{c|c}
-\RH_j & 0 \\ \hline
0 & -\RH_j
\end{array}
\right)\enspace,\qquad &
\sigma_{34}&=\left(
\begin{array}{c|c}
\RH_i & 0 \\ \hline
0 & \RH_i
\end{array}
\right)\enspace, \\
\end{aligned}
}
\end{equation}

\begin{equation*}\label{eq:Jspin52}
\scriptsize{\begin{aligned}
\sigma_{15}=\left(
\begin{array}{c|c}
0 & -\Id \\ \hline
\Id & 0
\end{array}
\right)\enspace,\quad
\sigma_{25}&=\left(
\begin{array}{c|c}
0 & \RH_i \\ \hline
\RH_i & 0
\end{array}
\right)\enspace,\quad &
\sigma_{35}&=\left(
\begin{array}{c|c}
0 & \RH_j \\ \hline
\RH_j & 0
\end{array}
\right)\enspace,\quad &
\sigma_{45}&=\left(
\begin{array}{c|c}
0 & \RH_k \\ \hline
\RH_k & 0
\end{array}
\right)\enspace.
\end{aligned}
}
\end{equation*}

\vspace{0.3cm}

\section{Structures defined on Grassmannians} In the previous Section we discussed how the groups $\Spin{8}, \Spin{6} \subset \SO{16}$ and $\Spin{5} \subset \SO{8}$ act on $\RR^{16}$ and $\RR^8$. This relates with the following definition of  Clifford system on a Riemannian manifold, cf. \cite[page 280]{ppv}.

\begin{de} A \emph{Clifford system $C_m$ on a Riemannian manifold $(M^N,g)$} is the datum of a rank $m+1$ vector sub-bundle $E^+ \subset \Endplus{TM}$, locally generated by a set of $m+1$ self-dual anti-commuting involutions that are related in the intersections of trivializing open sets by matrices in $\SO{m+1}$.
\end{de}

Here, for $\Spin{8}, \Spin{6} \subset \SO{16}$ and on $M^N=\RR^N$, the point of view of a Clifford system $C_m$, set of self-dual anti-commuting involutions, is equivalent to the approach we used in the previous Section, namely of anti-commuting complex structures $m_u$, see also Remark \ref{re} and Propositions \ref{8involutions}, \ref{6involutions}. The three cases $\Spin{8}, \Spin{6}$ and $\Spin{5}$ correspond to $m=7,5,4$ and $N=16,16,8$, cf. \cite[Table D]{ppv}.

This, together with the decomposition $T Gr \cong W \otimes W^\perp,$ with the holonomy of the three Grassmannians and with their spin property, gives the following:

\begin{lm}\label{lemma}
The Grassmannians 
\begin{equation}
Gr_{8}(\RR^{2n+8})=\frac{\mathrm{SO}(2n+8)}{\mathrm{SO}(2n)\times \mathrm{SO}(8)},  \quad Gr_{4}(\CC^{2n+4})=\frac{\mathrm{SU}(2n+4)}{\mathrm{S}(\mathrm{U}(2n)\times \mathrm{U}(4))}, \quad  Gr_{2}(\HH^{m+2})=\frac{\mathrm{Sp}(m+2)}{\mathrm{Sp}(m)\times \mathrm{Sp}(2)}
\end{equation}
admit a Clifford system respectively $C_7$, $C_5$, $C_4$ as in the previous definition. 

The first two Grassmannians admit also a rank $m+1$ vector subbundle $E^- \subset \Endminus{TM}$ locally generated by systems of anti-commuting complex structures $m_u$, where $u=1, \dots ,8$, $u=1, \dots , 6$, respectively, and whose compositions give rise to the same local almost complex structures as compositions of the corresponding local involutions.
\end{lm}

\begin{re}When, for the first two series of Grassmannians, the ambient space has odd dimension, i. e. for
\begin{equation}
Gr_{8}(\RR^{2n+9})=\frac{\mathrm{SO}(2n+9)}{\mathrm{SO}(2n+1)\times \mathrm{SO}(8)},  \quad Gr_{4}(\CC^{2n+5})=\frac{\mathrm{SU}(2n+5)}{\mathrm{S}(\mathrm{U}(2n+1)\times \mathrm{U}(4))}, 
\end{equation}
the manifold is not spin, cf. \cite{bh}, \cite{str}. As a consequence, one can still define the vector sub-bundles $\Endplus{T \, Gr}, \Endminus{T \, Gr}$ of the previous lemma, but only locally. This situation motivates the requiring of only a locally defined vector bundle $E$ in the definition of even Clifford structure, cf. the Introduction.
\end{re}

\vspace{0.3cm}

\section{The oriented Grassmannian $Gr_{8}(\RR^{2n+8})$}\label{realgrassmannians}

Let $n \geq 1$ and let $W$ be the rank $8$ tautological vector bundle over the Grassmannian $Gr_{8}(\RR^{2n+8})$ of oriented $8$-planes in $\RR^{2n+8}$. Let $W^\perp$ be the rank $2n$ orthogonal complement of $W$ in the trivial bundle given by the ambient linear space. Both $W$ and $W^\perp$ are Euclidean oriented vector bundles with the induced Euclidean metric. The tangent bundle of the Grassmannian decomposes as
\[
T\, Gr_{8}(\RR^{2n+8}) \cong \Hom{W,W^\perp} \cong W \otimes W^\perp.
\]

Thus, from local orthonormal bases $w_1 \dots , w_8$ and $w_1^\perp , w_2^\perp , \dots , w_{2n-1}^\perp, w_{2n}^\perp$
of sections respectively of $W$ and of $W^\perp$, one gets the following local basis of tangent vectors of $Gr_{8}(\RR^{2n+8})$:
\begin{equation}\label{thex}
\begin{aligned}
&x_{1,1}=w_1 \otimes w_1^\perp , &\dots \dots \;  & &\quad x_{8,1}= w_8 \otimes w_1^\perp ,\\
&x_{1,2}=w_1 \otimes w_2^\perp , &\dots \dots \;  & &\quad x_{8,2}= w_8 \otimes w_2^\perp ,\\
& \dots \dots \dots \dots \dots &\dots \dots \;  && \quad \dots \dots \dots \dots \dots \\
&x_{1,2n-1}=w_1 \otimes w_{2n-1}^\perp , &\dots \dots \; & &\quad x_{8,2n-1}= w_8 \otimes w_{2n-1}^\perp , \\
&x_{1,2n}=w_1\otimes w_{2n}^\perp  ,&\dots \dots \; & &\quad x_{8,2n}= w_8 \otimes w_{2n}^\perp .
\end{aligned}
\end{equation}

These 8-ples of sections can be written formally as octonions, i.e.  for all $\alpha=1, \dots , 2n$: 

\begin{equation}\label{thexx}
\vec x_\alpha= (x_{1,\alpha}, x_{2, \alpha }, \dots , x_{8,\alpha }) = x_{1,\alpha }+ix_{2,\alpha }+ \dots + hx_{8,\alpha},
\end{equation}
and can be ordered as a $n$-ple of pairs of octonions:
\[
\big( (\vec x_1, \vec x_2), \dots , (\vec x_{2n-1}, \vec x_{2n}) \big) \in (\OO \oplus \OO)^n.
\]

The even Clifford structure on $Gr_{8}(\RR^{2n+8})$ is defined by looking at a rank 8 Euclidean vector bundle $E \subset \Endminus {T\, Gr_{8}(\RR^{2n+8})}$, satisfying the conditions to be locally generated by anti-commuting orthogonal complex structures, that we denote here by $m_1, m_2, \dots , m_8$, in correspondence with the  $m_1,m_i,\dots,m_h$ of formulas \eqref{top}. The existence of such a $E$ is insured by the holonomy structure $\mathrm{SO}(2n)\times \mathrm{SO}(8)$ of the Grassmannian, by its spin property, and by the description of $\Spin{8}$ given in Section 2, cf. Lemma \ref{lemma}.

Accordingly, if $u,v$ are local sections of $E$, we can look at them as octonions in the basis $m_1, m_2, \dots , m_8$. For any such orthonormal pair $(u,v)$, look at $u \wedge v$ as a section of $\Lambda^2 E$, and define:
\[
\varphi: \Lambda^2 E \rightarrow \Endminus {T\, Gr_{8}(\RR^{2n+8})}
\] 
by
\begin{equation}\label{ourphi}
\varphi (u \wedge v) \big( (\vec x_1, \vec x_2), \dots , (\vec x_{2n-1}, \vec x_{2n}) \big) = \big( m_{u,v} (\vec x_1, \vec x_2), \dots , m_{u,v}( \vec x_{2n-1}, \vec x_{2n}) \big),
\end{equation}
i.e. by applying diagonally the matrix \eqref{muv}. 

Extending by Clifford composition this gives the Clifford morphism
\[
\varphi: \; \mathrm{Cl}^0 E \rightarrow \End{T\, Gr_{8}(\RR^{2n+8})}.
\] 

Thus, the 28 matrices listed in \eqref{eq:J1} describe how $\varphi$ associates, to the natural basis of the local sections of $\Lambda^2 E$, local almost complex structures on the Grassmannian $Gr_{8}(\RR^{2n+8})$. This set of 28 almost complex structures can be viewed as the local ingredient of the rank $8$ even Clifford structure. Note that this definition of $\varphi$ is coherent with that described in \cite[pages 950 and 953]{ms}, by looking at the action of the holonomy group ${\mathrm{SO}(2n)\times \mathrm{SO}(8)}$ of $Gr_{8}(\RR^{2n+8})$ on the model tangent space $\RR^{16n}$. Also, the orthogonal representation 
\[
\mathrm{SO}(2n)\times \mathrm{SO}(8) \rightarrow \mathrm{SO}(16n)
\]
defines an equivariant algebra morphism $\varphi: \mathrm{Cl}^0_8 \rightarrow \End{\RR^{16n}}$ mapping $\mathfrak{so}(8)=\mathfrak{spin}(8) \subset \mathrm{Cl}^0_8$ into $\frak{so}(16n) \subset \End{\RR^{16n}}$. This is insured by the triality of the $\mathfrak{spin}(8)$ representations, outer automorphism of $\Spin{8}$ interchanging its three non-equivalent representations in $\RR^8$, cf. Remark \ref{triality}. The parallel non-flat feature of the map $\varphi$
follows from the holonomy based construction.

This is summarized in the following:

\begin{te} There is a rank 8 vector sub-bundle $E \subset \Endminus {T\, Gr_{8}(\RR^{2n+8})}$ locally generated by anti-commuting orthogonal complex structures $m_1, m_2, \dots , m_8$, and $E$ defines on $Gr_{8}(\RR^{2n+8})$ an even non-flat parallel Clifford structure of rank 8. The morphism
\[
\varphi: \; \mathrm{Cl}^0 E \rightarrow \End{T\, Gr_{8}(\RR^{2n+8})}
\] 
is given by the Clifford extension of the map:  
\[
u \wedge v \in \Lambda^2 E \longrightarrow [m_{u,v}: (\OO \oplus \OO)^n \rightarrow (\OO \oplus \OO)^n],
\]
defined by applying diagonally the matrix \eqref{muv}. Here $u,v$ are local orthonormal sections of $E$, thus unitary orthogonal octonions in the basis $m_1, m_2, \dots , m_8$, so that $m_{u,v}$ acts diagonally on the $n$-ples of pairs of local tangent vectors:
\[
\big( (\vec x_1, \vec x_2), \dots , (\vec x_{2n-1}, \vec x_{2n}) \big),
\]
that according to formulas \eqref{thex}, \eqref{thexx} can be looked at as elements of $(\OO \oplus \OO)^n$.
\end{te}

\vspace{0.3cm}

\section{The complex Grassmannian $Gr_{4}(\CC^{2n+4})$}\label{complexgrassmannians}

Again with $n \geq 1$, let $W$ be here the rank $4$ tautological complex vector bundle over the complex Grassmannian $Gr_{4}(\CC^{2n+4})$, and let $W^\perp$ be its orthogonal complement in the trivial bundle, given by the linear space $\CC^{2n+4}$. The decomposition
\[
T\, Gr_{4}(\CC^{2n+4}) \cong \Hom{W,W^\perp} \cong W \otimes W^\perp,
\]
allows to write local frames, as follows.

Let $w_1, w_2 , w_3, w_4 $ and $w_1^\perp , w_2^\perp , \dots , w_{2n-1}^\perp, w_{2n}^\perp $
be local orthonormal bases of sections of $W$ and of $W^\perp$. Then a local orthonormal frame of tangent vectors on $Gr_{4}(\CC^{2n+4})$ is
\begin{equation}\label{thez}
\begin{aligned}
&z_{1,1}=w_1 \otimes w_1^\perp , &\dots \dots \;  & &\quad z_{4,1}= w_4 \otimes w_1^\perp ,\\
&z_{1,2}=w_1 \otimes w_2^\perp , &\dots \dots \;  & &\quad z_{4,2}= w_4 \otimes w_2^\perp ,\\
& \dots \dots \dots \dots \dots &\dots \dots \;  && \quad \dots \dots \dots \dots \dots \\
&z_{1,2n-1}=w_1 \otimes w_{2n-1}^\perp , &\dots \dots \; & &\quad z_{4,2n-1}= w_4 \otimes w_{2n-1}^\perp , \\
&z_{1,2n}=w_1\otimes w_{2n}^\perp  ,&\dots \dots \; & &\quad z_{4,2n}= w_4 \otimes w_{2n}^\perp .
\end{aligned}
\end{equation}

Again, one can look at the above lines as ($\alpha=1, \dots , 2n$) 
\begin{equation}\label{thezz}
\vec z_\alpha= (z_{1,\alpha}, z_{2, \alpha }, z_{3, \alpha} z_{4,\alpha }) \in \CC^4,
\end{equation}
and order them as an $n$-ple of pairs:
\[
\big( (\vec z_1, \vec z_2), \dots , (\vec z_{2n-1}, \vec z_{2n}) \big) \in (\CC^4 \oplus \CC^4)^n.
\]

Look now at the vector sub-space $F=<1,i,j,k,e,f> \subset \OO$, and recall that the corresponding operators $m_u$, with $u \in F$, act on the complex vector space $\CC^4$, cf. formulas \eqref{complex}.

Similarly to what described in the previous Section for the real Grassmannian, there is a vector sub-bundle $E \subset \Endminus{T \; Gr_4(\CC^{2n+4})}$, locally generated by anti-commuting orthogonal complex structures $m_1,m_2, \dots , m_6$, and corresponding to the $m_1,m_i,m_j,m_k,m_e,m_f$ of Section 2. This is due to the holonomy $\mathrm{S(U}(2n)\times \mathrm{U}(4))$ of the Grassmannian and its spin property, that together allow the description of $\mathrm{SU}(4)\cong \Spin{6}$ given in Section 2, cf \ref{lemma}. If $(u,v)$ is an orthonormal pair of sections of $E$, then $u \wedge v$ is a section of $\Lambda^2 E$, and the map
\[
\varphi: \Lambda^2 E \rightarrow \Endminus {T\, Gr_{4}(\CC^{2n+4})},
\] 
given by
\begin{equation}\label{ourphi}
\varphi (u \wedge v) \big( (\vec z_1, \vec z_2), \dots , (\vec z_{2n-1}, \vec z_{2n}) \big) = \big( m_{u,v} (\vec z_1, \vec z_2), \dots , m_{u,v}( \vec z_{2n-1}, \vec z_{2n}) \big),
\end{equation}
extends by Clifford composition to the Clifford morphism
\[
\varphi: \; \mathrm{Cl}^0 E \rightarrow \End{T\, Gr_{4}(\CC^{2n+4})}.
\] 

Thus, among the matrices listed in \eqref{eq:J1}, to describe $\varphi$ we take now into account only the 15 local almost complex structures on $Gr_{4}(\CC^{2n+4})$ given by choices $u,v \in \{1,i,j,k,e,f\}$. This set of 15 almost complex structures gives the local rank $6$ even Clifford structure, as considered in \cite[pages 950 and 953]{ms}. Note that the holonomy group ${\mathrm{S(U}(2n)\times \mathrm{U}(4))}$ acts on the model tangent space $\CC^{8n}$, and the orthogonal representation 
\[
\mathrm{S(U}(2n)\times \mathrm{U}(4)) \rightarrow \mathrm{SU}(8n)
\]
defines an equivariant algebra morphism $\varphi: \mathrm{Cl}^0_6 \rightarrow \End{\CC^{8n}}$ mapping $\mathfrak{su}(4)=\mathfrak{spin}(6) \subset \mathrm{Cl}^0_6$ into $\frak{su}(8n) \subset \End{\CC^{8n}}$. The parallel non-flat feature of $\varphi$
follows again from the holonomy based construction. This gives the following:

\begin{te} There is a rank 6 vector sub-bundle $E \subset \Endminus {T\, Gr_{4}(\CC^{2n+4})}$ locally generated by anti-commuting orthogonal complex structures $m_1, m_2, \dots , m_6$, and $E$ defines on $Gr_{4}(\CC^{2n+4})$ an even non-flat parallel Clifford structure of rank 6. The morphism
\[
\varphi: \; \mathrm{Cl}^0 E \rightarrow \End{T\, Gr_{4}(\CC^{2n+4})}
\] 
is given by Clifford extension of the map:  
\[
u \wedge v \in \Lambda^2 E \longrightarrow [m_{u,v}: (\CC^4 \oplus \CC^4)^n \rightarrow (\CC^4 \oplus \CC^4)^n],
\]
defined by diagonal extension of the matrix \eqref{muv}. Here $u,v$ are local orthonormal sections of $E$, thus unitary orthogonal in the basis $m_1, m_2, \dots , m_6$, so that $m_{u,v}$ acts diagonally on the $n$-ples of pairs of local tangent vectors:
\[
\big( (\vec z_1, \vec z_2), \dots , (\vec z_{2n-1}, \vec z_{2n}) \big),
\]
that according to formulas \eqref{thez}, \eqref{thezz} can be looked at as elements of $(\CC^4 \oplus \CC^4)^n$.
\end{te}

\vspace{0.3cm}

\section{The quaternionic Grassmannian $Gr_{2}(\HH^{n+2})$}\label{quaternionicgrassmannians}

Here the vanishing of the integral cohomology $H^2$ insures that of the Stiefel-Whitney class $w_2$ and for any $n$ the spin property of the Grassmannian. By the decomposition
\[
T\, Gr_{2}(\HH^{n+2}) \cong \Hom{W,W^\perp} \cong W \otimes W^\perp,
\]
from orthonormal local frames $w_1, w_2$ of $W$ and $w_1^\perp , w_2^\perp , \dots , w_{n-1}^\perp, w_{n}^\perp $ of $W^\perp$
one gets the local orthonormal frame of tangent vectors
\begin{equation}\label{theh}
\begin{aligned}
&h_{1,1}=w_1 \otimes w_1^\perp ,\quad &h_{2,1}= w_2 \otimes w_1^\perp ,&\\
&h_{1,2}=w_1 \otimes w_2^\perp ,\quad &h_{2,2}= w_2 \otimes w_2^\perp ,&\\
& \dots \dots \dots \dots \dots  \quad &\dots \dots \dots \dots \dots &\\
&h_{1,n-1}=w_1 \otimes w_{n-1}^\perp ,\quad &h_{2,n-1}= w_2 \otimes w_{n-1}^\perp , &\\
&h_{1,n}=w_1\otimes w_{n}^\perp  ,\quad &h_{2,n}= w_2\otimes w_{n}^\perp . &
\end{aligned}
\end{equation}
Then write ($\alpha=1, \dots , n$) 
\begin{equation}\label{thehh}
\vec h_\alpha= (h_{1,\alpha}, h_{2, \alpha }) \in \HH \oplus \HH,
\end{equation}
and order the $\vec h_\alpha$ as a $n$-ple
\[
\big( \vec h_1, \vec h_2, \dots \vec h_{n-1}, \vec h_{n} \big) \in (\HH \oplus \HH)^{n}.
\]

Thus, taking into account the $\Spin{5}$ subsection of Section \ref{Spin8}:

\begin{te} There is a rank 5 vector sub-bundle $E \subset \Endplus {T\, Gr_{2}(\HH^{n+2})}$ locally generated by anti-commuting orthogonal self-dual involutions $\sigma_1, \sigma_2, \dots , \sigma_5$, and $E$ gives rise to an even non-flat parallel Clifford structure of rank 5  on $Gr_{2}(\HH^{n+2})$. The Clifford morphism $\varphi$ is constructed as follows. Let $u,v$ be local orthonormal sections of $E$, and let the composition $uv$ act diagonally on the $n$-ples:
\[
\big( \vec h_1, \vec h_2, \dots , \vec h_{n-1}, \vec h_{n} \big),
\]
that according to formulas \eqref{theh}, \eqref{thehh} can be looked at as elements of $(\HH \oplus \HH)^n$.
Then the morphism
\[
\varphi: \; \mathrm{Cl}^0 E \rightarrow \End{T\, Gr_{2}(\HH^{n+2})}
\] 
is given by the Clifford extension of the map:  
\[
u \wedge v \in \Lambda^2 E \longrightarrow [uv: (\HH \oplus \HH)^{n} \rightarrow  (\HH \oplus \HH)^{n} ].
\]
\end{te}

\vspace{0.3cm}

\section{$16$-dimensional examples}

\subsection{The complex Cayley line $Gr_8(\RR^{10})$} \label{firstgrassmannian} 

The complex K\"ahler structure of $Gr_8 (\RR^{10})$ can be recognized via the diffeomorphism with the non singular quadric $Q_8 \subset \CC P^9$, defined by choosing a complementary oriented 2-plane $\vec p \wedge \vec q$, fibre of $W^\perp$ ($\vec p =(p_0, \dots , p_9), \vec q =(q_0,\dots q_9) \in \RR^{10}$ orthonormal), so that $\vec z = \vec p +i \vec q$ satisfies the equation $z_0^2+ \dots + z_9^2$ of $Q_8$. Locally, the complex structure $I$ of $Gr_8 (\RR^{10})$ is given, in the notations of formula
\eqref{thex}, by 
\[
I: x_{1,1} \rightarrow x_{2,1},\; I: x_{1,2} \rightarrow x_{2,2}, \qquad \dots \dots \qquad  I: x_{1,8} \rightarrow x_{2,8}.
\]
The following Proposition describes families of sub-manifolds that, according to the classification in \cite{cn}, are maximal totally geodesic in $Gr_8 (\RR^{10})$.

\begin{pr} The Grassmannian $Gr_8 (\RR^{10}) \cong Q_8 \subset \CC P^9$, contains two families $\mathcal F, \mathcal F'$ of complex projective $4$-spaces, respectively $\CC P^4, \CC {P'^4}$, each family being parametrized by the Hermitian symmetric space $\SO{10}/\U{5}$ of the linear complex structures on $\RR^{10}$. There are also two families of $8$-spheres $\mathcal H =\{\OO P^1 \cong S^8\}$, $\mathcal H' =\{\OO P'^1 \cong S'^8\}$, both parametrized by the sphere $S^9$. Each point of $Gr_2 (\RR^{10})$ belongs to infinitely many projective subspaces $\CC P^4, \CC {P'^4}$  and infinitely many spheres $S^8, S'^8$ of the families $\mathcal F$, $\mathcal F'$, $\mathcal H$, $\mathcal H'$.
\end{pr}
\bigskip

The first two families are in fact in correspondence with the linear complex structures $J: \RR^{10} \rightarrow \RR^{10}$: for each $J$ there is a family of $2$-planes in $\RR^{10}$ that are invariant under $J$, i.e. a family of complex lines through the origin in the $5$-dimensional complex vector space $(\RR^{10}, J)$. There are two families of such complex structures $J$, according the orientation they induce on $\RR^{10}$, and each of the two families is parametrized by the Hermitian symmetric space $\SO{10}/\U{5}$. Indeed the full set of "orthogonal complex structures" on $\RR^{10}$ is the homogeneous space $\mathrm{O}(10)/\U{5}$, consisting of two copies of the mentioned Hermitian symmetric space, cf. \cite{Sa}. Concerning the families of 8-spheres, they correspond to the oriented lines in the possible oriented $\RR^9 \subset \RR^{10}$ and to the oriented 2-planes that include a fixed oriented line in $\RR^{10}$.

The real cohomology of $Gr_8 (\RR^{10})$ is fully in line with its two parallel non-flat Clifford structures: the one of rank 2, giving the complex K\"ahler structure, and the rank 8, described in Section \ref{realgrassmannians}:

\begin{equation}\label{rosenfeld1}
H^*(Gr_8 (\RR^{10});\RR) \cong \RR[e, e^\perp]/(\rho_{5},\rho_{8})
\end{equation}
where $e \in H^8$ and $e^\perp \in H^2$ are the Euler classes of the vector bundles $W$ and of $W^\perp$, and the relations are $\rho_5=ee^\perp \in H^{10}, \; \rho_8=e^2 - (e^\perp)^8 \in H^{16}$. 
Its Poincar\'e polynomial reads:
\[
\mathrm{Poin}_{_{Gr_2 (\RR^{10})}} = 1+t^2+t^4+t^6+2t^8+t^{10}+t^{12}+t^{14}+t^{16} ,
\]
and as representatives of the two generators of the cohomology algebra one can choose the complex K\"ahler form $\omega$ and an "octonionic K\"ahler" $8$-form $\Psi$, that comes from the Chern-Weil representation of $e$.

In the octonionic setting, it is natural to look at $Gr_8 (\RR^{10})$ as the complex octonionic projective line $(\CC \otimes \OO)P^1$ and as such a totally geodesic half- dimensional sub-manifold of $\EIII$, the complex octonionic projective plane. Then the 8-form $\Psi$ is related to a construction involving the holonomy group $\Spin{10} \cdot \U{1}$ of $\EIII$, cf. \cite{pp3}.

\vspace{0.2cm}

\subsection{The third Severi variety $Gr_4(\CC^6)$}

This is a remarkable Grassmannian, supporting three different non-flat parallel even Clifford structures, of rank $2,3,6$ respectively. The first two of them come from the complex K\"ahler and the from the quaternion K\"ahler structure of $Gr_4(\CC^6) \cong Gr_2(\CC^6)$, and the third one comes from the construction in Section \ref{complexgrassmannians}. The name we use here has been proposed by F. Zak \cite{ZakSeV}, who showed that $Gr_4(\CC^6)$, sub-variety of $\CC P^{14}$ in its Pl\"ucker embedding, is one of the four smooth projective algebraic sub-varieties of critical codimension in a $\CC P^{N}$, not contained in a hyperplane, and unable to fill the ambient projective ambient space through their secant and tangent lines.

The real cohomology is generated by the complex K\"ahler and the quaternion K\"ahler 4-form, and:
\[
\mathrm{Poin}_{_{Gr_4 (\CC^{6})}} = 1+t^2+2t^4+2t^6+3t^8+2t^{10}+2t^{12}+t^{14}+t^{16} .
\]

Like $\EIII$, also $Gr_4(\CC^6)$ appears in a table of projective planes over composition algebras, cf. \cite{pp3}. As such $Gr_4(\CC^6) \cong (\CC \otimes \HH)P^2$, and its corresponding totally geodesic "complex quaternionic" projective line is $ (\CC \otimes \HH)P^1 \cong Gr_2(\CC^4) \cong Gr_4(\RR^6) \cong Q_4 \subset \CC P^5$.

\vspace{0.2cm}

\subsection{The  quaternionic Klein quadric $Gr_2(\HH^4)$}

This Grassmannian, that can be also viewed as parametrizing all projective lines $l$ in a $\HH P^3$, is neither K\"ahler nor quaternion K\"ahler. Although the point of view of the Pl\"ucker coordinates over the skew-field $\HH$ is not straightforward and it is not clear how to represent $Gr_2(\HH^4)$ in an ambient quaternionic projective space, nevertheless $Gr_2(\HH^4)$ shares the following features of the classical non singular Klein quadric $Q_4 \subset P^5$. Cf. also the similar properties in this Section for the higher dimensional complex quadric $Q_8 \cong Gr_8 (\RR^{10})$.

\begin{pr} $Gr_2(\HH^4)$ contains two families $\mathcal F, \mathcal F'$ of $2$-planes, respectively $\HH P^2, \HH {P'^2}$, each family being parametrized by a $\HH P^3$. The planes $\HH P^2$ of the family $\mathcal F$ represent the "ruled planes" in $\HH P^3$, i. e. the sets of lines $l$ that are contained in a projective plane of $\HH P^3$. The planes $\HH {P'^2}$ of the family $\mathcal F'$ represent the "stars of lines" in $\HH P^3$, i. e. the sets of lines $l$ that contain a point of $\HH P^3$. It follows that planes of the same family $\mathcal F$ or $\mathcal F'$ intersect in a point of $Gr_2(\HH^4)$, and that planes of different families either do not intersect of intersect in a line contained in $Gr_2(\HH^4)$. 
\end{pr}

The Poincar\'e polynomial:
\[
\mathrm{Poin}_{_{Gr_2 (\HH^{4})}} = 1+t^4+2t^8+t^{12}+t^{16},
\]
shows as generators of the cohomology the symplectic Pontrjagin classes of the quaternionic vector bundles $W$ and $W^\perp$, cf \cite{p0}.

\vspace{0.3cm}

\section{Some other examples}

\subsection{The Wolf space $Gr_8(\RR^{12})$}

This Grassmannian has also three non-flat even Clifford structures. There are two of rank 3, corresponding to the two quaternion K\"ahler structures (in correspondence with two different ways to define hypercomplex structures on the planes on any $Gr_4(\RR^{n+4})$), and the one of rank 8, described in Section \ref{realgrassmannians}. 
Indeed, $Gr_8(\RR^{12})$ can be looked at as the "quaternion octonionic" projective line $(\HH \otimes \OO)P^1$, totally geodesic submanifold of the exceptional symmetric space $(\HH \otimes \OO)P^2 \cong \EVI$, cf. \cite{e}. 

Its Poincar\'e polynomial
\[
\mathrm{Poin}_{_{Gr_4 (\RR^{12})}} = 1+2t^4+4t^8+5t^{12}+6t^{16}+5t^{20}+4t^{24}+2t^{28}+t^{32}
\]
exhibits the presence of two quaternion K\"ahler 4-forms and of an "octonionic K\"ahler" $8$-form $\Psi$. This latter $\Psi$ is related to one that is defined on $\EVI$ through its holonomy group $\Spin{12} \cdot \Sp{1}$, cf. \cite{p}.

\vspace{0.2cm}

\subsection{The two 64-dimensional examples $Gr_8(\RR^{16})$ and $Gr^\perp_8(\RR^{16})$}

On the Grassmannian $$Gr_8(\RR^{16})= \frac{\SO{16}}{\SO{8} \times \SO{8}}$$ one can define, besides the parallel even Clifford structure of rank 8 described in Section \ref{realgrassmannians}, another similar structure obtained by interchanging the r\^ole of the two vector bundles $W$ and $W^\perp$, i. e. by operating through the $m_{u,v}$ on elements of $W^\perp$. The real cohomology
\begin{equation}\label{coho}
H^*(Gr_8(\RR^{16})) \cong \frac{\mathbb R [e, p_1, p_2, p_3, e^\perp, p_1^\perp, p_2^\perp, p_3^\perp]}{ee^\perp =0, \; (1+p_1+ p_2+ p_3)(1+p_1^\perp+ p_2^\perp,+p_3^\perp)=1} \; ,
\end{equation}
in terms of Pontrjagin classes $p_\alpha, p_\alpha^\perp$ and Euler classes $e, e^\perp$ of $W$ and $W^\perp$, gives rise to the Poincar\'e polynomial
\[
\mathrm{Poin}_{_{Gr_8 (\RR^{16})}} = 1+t^4+4t^8+5t^{12}+9t^{16}+11t^{20}+15t^{24}+15t^{28}+18t^{32}+ \dots \; ,
\]
and hence to the Euler characteristic $\chi (Gr_8(\RR^{16})= 140$, cf. \cite[Table I, page 339]{hs}.

Due to their construction, the two mentioned even Clifford structures descend to a unique even parallel Clifford structure of rank 8 on the smooth $\mathbb Z_2$-quotient
\[
Gr^\perp_8(\RR^{16}) = Gr_8(\RR^{16})/\perp
\]
by the orthogonal complement involution $\perp$. This quotient $Gr^\perp_8(\RR^{16})$ turns out to be a totally geodesic half dimensional submanifold of $\EVIII$, the largest of the exceptional symmetric spaces of compact type, and can be looked at as the "projective line" $(\OO \otimes \OO)P^1$ over the "octonionic octonions", cf. \cite{cn}, \cite{e}.

For the computation of the cohomology of $Gr^\perp_8(\RR^{16})$, just note that the involution $\perp$ identifies $p_1 \rightarrow p_1^\perp, p_2 \rightarrow p_2^\perp, p_3 \rightarrow p_3^\perp, e \rightarrow e^\perp$, cf. the similar computation in \cite{p0.5} for the cohomology of the quaternion K\"ahler analogue $Gr^\perp_4(\RR^{8})$ . This, due to the relations in \eqref{coho}, allows to survive, up to dimension 32, only the classes $p_1^2,e,p_1^4, p_1^2 e,p_1^6, p_1^4 e,p_1^8,p_16 e$. This gives now the Poincar\'e polynomial
\[
\mathrm{Poin}_{_{Gr^\perp_8 (\RR^{16})}} = 1+2t^8+2t^{16}+2t^{24}+2t^{32}+ \dots \; ,
\]
and the Euler characteristic $\chi(Gr^\perp_8 (\RR^{16}) = 16.$

\end{document}